\documentclass[11pt,noinfoline]{imsart}
\usepackage{amsmath,amsfonts,amsthm}
\usepackage{graphicx}

\newtheorem{thm}{Theorem}[section]
\newtheorem{theorem}{Theorem}[section]
\newtheorem{lem}{Lemma}[section]

\newtheorem{proposition}[thm]{Proposition}
\newtheorem{cor}[thm]{Corollary}
\newtheorem{prop}[thm]{Proposition}

\newtheorem{rem}{Remark}[section]
\newtheorem{example}{Example}[section]

\begin{document}

\def\l{\lambda}
\def\t{\theta}
\def\T{\Theta}
\def\m{\mu}
\def\a{\alpha}
\def\b{\beta}
\def\g{\gamma}
\def\o{\omega}
\def\p{\varphi}
\def\D{\Delta}
\def\O{\Omega}
\def\G{\Gamma}

\def\N{{\mathbb N}}
\def\C{{\mathbb C}}
\def\Z{{\mathbb Z}}
\def\R{{\mathbb R}}
\def\P{{\mathbb P}}
\def\E{{\mathbb E}}

\def\L{{\mathcal L}}
\def\D{{\mathcal D}}
\def\I{{\mathcal I}}
\def\K{{\mathbb K}}
\def\T{{\mathcal T}}

\begin{frontmatter}

\title{Exponential functionals of Brownian motion and class-one Whittaker functions}
\runtitle{Exponential functionals and Whittaker functions}

\begin{aug}

\author{\fnms{Fabrice}  \snm{Baudoin}\corref{}
  \ead[label=e1]{fbaudoin@math.purdue.edu}}

  \affiliation{Purdue University}
  
   \address{Department of Mathematics, Purdue University, West Lafayette, IN 47906, USA\\ 
          \printead{e1}}
          
  \author{\fnms{Neil}  \snm{O'Connell}\corref{}
  \ead[label=e2]{n.m.o-connell@warwick.ac.uk}}

  \runauthor{F. Baudoin and N. O'Connell}

  \affiliation{University of Warwick}

  \address{Mathematics Institute, University of Warwick, Coventry CV4 7AL, UK\\ 
          \printead{e2}}

\end{aug}

\begin{abstract}
We consider exponential functionals of a Brownian motion with drift in $\R^n$, 
defined via a collection of linear functionals.  We give a characterisation of the Laplace
transform of their joint law as the unique bounded solution, up to a constant factor, 
to a Schr\"odinger-type partial differential 
equation. We derive a similar equation for the probability density.  
We then characterise all diffusions which can be interpreted as having the law of the
Brownian motion with drift conditioned on the law of its exponential functionals.  In the
case where the family of linear functionals is a set of simple roots, the Laplace transform 
of the joint law of the corresponding exponential functionals can be expressed in terms of 
a (class-one) Whittaker function associated with the corresponding root system. 
In this setting, we establish some basic properties of the corresponding diffusion processes.
\end{abstract}



\end{frontmatter}

\section{Introduction}

Let $(X_t,t\ge 0)$ be a standard one-dimensional Brownian motion with drift $\mu$.
In the paper~\cite{my}, Matsumoto and Yor consider the process 
\[\left( \log \int_0^t e^{2X_s-X_t}ds ,\ t > 0\right),\]
and prove that it is a diffusion process with infinitesimal generator given by 
\begin{equation}\label{myg}
\frac12 \frac{d^2}{dx^2}+\left(\frac{d}{dx}\log K_\mu(e^{-x})\right) \frac{d}{dx},
\end{equation}
where $K_\mu$ is the Macdonald function.  As explained in~\cite{my}, this theorem can
be regarded, by Brownian scaling and Laplace's method, as a generalization of Pitman's
`$2M-X$' theorem~\cite{pitman,rp} which states that, if $M_t=\max_{0\le s\le t} X_s$, then
$2M-X$ is a diffusion process with infinitesimal generator given by
\begin{equation}\label{rpg}
\frac12 \frac{d^2}{dx^2}+ \mu\coth(\mu x) \frac{d}{dx}.
\end{equation}

Note that $\mu\coth(\mu x)=(-\mu)\coth(-\mu x)$ and $K_\mu=K_{-\mu}$.  Suppose $\mu>0$.
Then the diffusion with generator (\ref{rpg}) can be interpreted as a Brownian motion 
with drift $\mu$ conditioned to stay positive.  Similarly, the diffusion with generator 
(\ref{myg}) can be interpreted as the Brownian motion $X$ conditioned on its 
exponential functional $A=\int_0^\infty e^{-2X_s} ds$
having a certain distribution (a Generalised Inverse Gaussian law) in a sense
which can be made precise~\cite{b,b2}.
The relevance of these interpretations in the present context is as follows.

Set $J=-\min_{t\ge 0} X_t$ and let $\tilde J$ be an independent
copy of $J$ which is also independent of $X$.  Then the process 
$\tilde X = 2 \max\{ M - \tilde J, 0 \} - X$ has the same law as $X$ and,
moreover, $\tilde J=-\min_{t\ge 0} \tilde X_t$.  This is well-known and can be seen
for example as a consequence of the classical output theorem for the $M/M/1$ 
queue~\cite{oy2}.  From this we can see that the process $2M-X$ has the same law 
as that of $X$ conditioned (in an appropriate sense) on the event that $J=0$; in
other words, $2M-X$ is a diffusion with infinitesimal generator given by (\ref{rpg}),
started from zero.  This basic idea can be used to obtain a multi-dimensional
version of Pitman's $2M-X$ theorem~\cite{bbo,bj,oy}, which gives a representation
of a Brownian motion conditioned to stay in a Weyl chamber in $\R^n$ as a certain 
functional (which generalizes $2M-X$) of a Brownian motion in $\R^n$.

Similarly~\cite{my}, if $\tilde A$ is an independent 
copy of $A$, which is also independent of $X$, then the process
\[\tilde X_t = \log\left(1+ \tilde A^{-1} \int_0^t e^{2(X_s-X_t)}ds\right) +X_t  ,\qquad t \ge 0\]
has the same law as $X$ and, moreover, $\tilde A = \int_0^\infty e^{-2\tilde X_s} ds$.
This time, we conclude that, for each $\epsilon>0$, the process
\[\log\left(\epsilon + \int_0^t e^{2(X_s-X_t)}ds\right) +X_t  ,\qquad t \ge 0\]
the same law as that of $\log\epsilon+X$ conditioned on the event $A=\epsilon$.
Carefully letting $\epsilon\to 0$ yields the theorem of Matsumoto and Yor~\cite{b}.
As explained in~\cite{oc03}, the probabilistic proofs of the multi-dimensional versions 
of Pitman's $2M-X$ theorem given in the papers~\cite{bbo,oy} carry over, in the same way, 
to the exponential functionals setting, although the task of letting the analogue of $\epsilon$ 
go to zero is a highly non-trivial problem in the general setting.  Nevertheless, it gives a 
heuristic derivation that a certain functional of a Brownian motion in $\R^n$ should have 
the same law as a Brownian motion conditioned on a certain collection of its exponential 
functionals.  This leads us to the question considered in the present paper.  

We consider a Brownian motion $B^{(\mu)}$ in $\R^n$ with drift $\mu$, and a 
collection of linear functionals $\a_1,\ldots,\a_d$ such that the exponential functionals
\[
A^i_\infty= \int_0^\infty e^{-2\a_i(B_s^{(\mu)})}ds \qquad (i=1,\ldots, d)
\]
are almost surely finite.  Our aim is to understand which diffusion processes can arise when
we condition on the law of $A_\infty=(A^1_\infty,\ldots,A^d_\infty)$.  The first step is to understand the law of $A_\infty$.
We show that the Laplace transform of $A_\infty$ satisfies a certain Schr\"odinger-type
partial differential equation and proceed to characterise all diffusion processes which 
can be interpreted as having the law of $B^{(\mu)}$ conditioned on the law 
of $A_\infty$. 

 In the case when
$\a_1,\ldots,\a_d$ is a simple system (see section 4 below for a definition), 
these diffusion processes are closely related to the 
quantum Toda lattice.  The Schr\"odinger operator is
\[H=\frac12 \Delta + \sum_i \theta_i^2 e^{-2\alpha_i},\]
where $\theta\in\R^d$, and the corresponding diffusion process has 
infinitesimal generator given by
\begin{equation}\label{blah}
\frac12 \Delta + \nabla \log k_\mu \cdot\nabla ,
\end{equation}
where $k_\mu$ is a particular eigenfunction of $H$ known as a class-one 
Whittaker function.  In the case $n=d=1$ and $\alpha_1(x)=x$, the class-one
Whittaker function is $k_\mu(x)=K_\mu(e^{-x})$ 
and the infinitesimal generator is given by (\ref{myg}).  
More generally, for a simple system $\a_1,\ldots,\a_d$, the diffusion process with
generator given by (\ref{blah}) plays an analogous role, in the exponential functionals setting, 
as that of a Brownian motion conditioned to stay in the Weyl chamber
$\{x\in\R^n:\ \a_i(x)>0,\ i=1,\ldots, d\}$. 
These processes have already found an application in the paper~\cite{noc10}, 
where the corresponding 
multi-dimensional version of the above theorem of Matsumoto and Yor in the `type A' case 
has been proved and used to determine the law of the partition function associated with a 
directed polymer model which was introduced in the paper~\cite{oy2}. 

The outline of the paper is as follows.  In section 2 we work in a general setting and establish a Schr\"odinger type partial differential equation satisfied by the characteristic function of exponential functionals of a multidimensional Brownian motion. We also study a family of martingales related to the conditional laws of exponential functionals that will later appear. In section 3, we identify a family
of diffusions which can be interpreted as having the law of the Brownian motion with drift conditioned 
on the law of its exponential functionals. In section 4, we restrict our attention to the case where the
collection of vectors used to define the exponential functionals is a simple system, and give an overview
of relevant facts about class-one Whittaker functions.  In section 5, we study properties of the 
conditioned processes in this setting.  In the final section, we present some explicit results for
the `type $A_2$' case.

\section{Exponential functionals and associated partial differential equations}

In this section, we work in a general setting and establish a Schr\"odinger type partial differential equation satisfied by the characteristic function of exponential functionals of a multidimensional Brownian motion. We also study a family of martingales related to the conditional laws of exponential functionals that will later appear.
 
Let $\a_1,\ldots,\a_d$ be a collection of distinct, non-zero vectors in $\R^n$
such that 
\begin{equation}
\Omega=\{x\in\R^n:\ \a_i(x)>0\ \forall i\}
\end{equation} 
is non-empty.
Let $B^{(\mu)}$ be a standard Brownian motion in $\R^n$ with drift $\mu\in\Omega$.
For $0\le t\le \infty$, set
\[
A^i_t= \int_0^t e^{-2\a_i(B_s^{(\mu)})}ds \qquad i=1,\ldots, d.
\]
Here, $\alpha_i(\beta)=(\alpha_i,\beta)$ where $(\cdot,\cdot)$ denotes
the usual inner product on $\R^n$.

\subsection{Partial differential equation for the characteristic function}

The process $\left( B_t^{(\mu)} ,A_t \right)_{t \ge 0}$ is a diffusion with generator
\[
\frac{1}{2} \Delta_x +(\mu,\nabla_x)+\sum_{i=1}^d e^{-2\alpha_i (x)} \frac{\partial }{\partial a_i }.
\]
We first check that this operator is hypoelliptic.
\begin{proposition}
The operator
\[
\frac{1}{2} \Delta_x +(\mu,\nabla_x)+\sum_{i=1}^d e^{-2\alpha_i (x)} \frac{\partial }{\partial a_i }
\]
is hypoelliptic on $\mathbb{R}^{n+d}$ and therefore, for $t>0$ the random variable $\left( B_t^{(\mu)} ,A_t \right)$ admits a smooth density with respect to the Lebesgue measure.
\end{proposition}
\begin{proof}
We use H\"ormander's theorem. Since the $\alpha_i$'s are pairwise different and non zero, there exists $v \in \mathbb{R}^n$ such that
\[
i \neq j \Rightarrow \alpha_i (v) \neq \alpha_j (v).
\]
Consider now the vector field
\[
V= \sum_{i=1}^n v_i \frac{\partial}{\partial x_i},
\]
and let us denote
\[
T=\sum_{i=1}^d e^{-2\alpha_i (x)} \frac{\partial }{\partial a_i }
\]
The Lie bracket between $V$ and $T$ is given by
\[
\mathcal{L}_V T=[V,T]=-2\sum_{i=1}^d \alpha_i (v) e^{-2\alpha_i (x)} \frac{\partial }{\partial a_i }.
\]
Similarly, by iterating this bracket $k$ times, we get
\[
\mathcal{L}^k_V T=(-1)^k 2^k\sum_{i=1}^d \alpha_i (v)^k e^{-2\alpha_i (x)} \frac{\partial }{\partial a_i }.
\]
Since the $\alpha_i$'s are pairwise different and non zero, we deduce from the Van der Monde determinant that at every $x \in \mathbb{R}^n$ the family
\[
\left\{ \mathcal{L}^k_V T, 1 \le k \le d  \right\}.
\]
is a basis of $\mathbb{R}^d$. It implies that the Lie bracket generating condition of H\"ormander is sastisfied so that the operator $\frac{1}{2} \Delta_x +(\mu,\nabla_x)+\sum_{i=1}^d e^{-2\alpha_i (x)} \frac{\partial }{\partial a_i }$ is hypoelliptic.
\end{proof}

Let now $\theta\in\R^d$ and, for $x\in \mathbb{R}^n$, define
\[
g_\mu^\theta (t,x) =\mathbb{E} \left( e^{-\sum_{i=1 }^d \theta_i^2 e^{-2\alpha_i (x) } A_t^i}\right), \quad t \ge 0,
\]
and
\[
j_\mu^\theta (x) =\mathbb{E} \left( e^{-\sum_{i=1 }^d \theta_i^2 e^{-2\alpha_i (x) } A_\infty^i}\right).
\]
\begin{proposition}
\begin{enumerate}
\item The semigroup generated by the Schr\"odinger operator
\[
\frac{1}{2} \Delta +(\mu,\nabla )-
 \sum_{i=1}^d \theta_i^2 e^{-2\alpha_i (x)}
\]
admits a heat kernel $q_\mu^\theta (t,x,y)$ and we have
\[
g_\mu^\theta (t,x)=\int_{\mathbb{R}^n} q_\mu^\theta (t,x,y) dy.
\]
\item
The function $j_\mu^\theta$ is the unique bounded function that satisfies the partial differential equation
\[
\frac{1}{2} \Delta j_\mu^\theta (x)+(\mu,\nabla j_\mu^\theta (x))=
\left( \sum_{i=1}^d \theta_i^2 e^{-2\alpha_i (x)}\right)j_\mu^\theta (x)
\]
and the limit condition
\[
\lim_{x\rightarrow \infty,x \in \Omega } j_\mu^\theta (x)=1.
\]
\end{enumerate}
\end{proposition}

\begin{proof}

\

\begin{enumerate}
\item
It is a straightforward consequence of the Feynman-Kac formula that $q_\mu^\theta (t,x,y)$ exists and is given by
\[
q_\mu^\theta (t,x,y)=\mathbb{E} \left( e^{-\sum_{i=1 }^d \theta_i^2 e^{-2\alpha_i (x) } A_t^i} \mid B_t^{(\mu)} =y-x \right) \frac{1}{(2\pi t)^{\frac{n}{2}}} e^{-\frac{\|y-x -\mu t\|^2}{2t}}.
\]
Integrating this with respect to $y$, we obtain
\[
g_\mu^\theta (t,x)=\int_{\mathbb{R}^n} q_\mu^\theta (t,x,y) dy.
\]
\item
It is again a straightforward consequence of the Feynman-Kac formula that $j_\mu^\theta$ 
solves the partial differential equation, and the limit condition is easily checked. 
Let us now prove uniqueness. 
We have to show that if $\phi$ is a bounded solution of the equation that satisfies
\[
\lim_{x\rightarrow \infty,x \in \Omega } \phi (x)=0,
\]
then $\phi=0$. For that, let us observe that under the above conditions, 
for $x\in \mathbb{R}^n$, the process
\[
\phi(B_t^{(\mu)}+x) \exp\left( -\sum_{i=1}^d \theta_i^2 e^{-2\alpha_i (x) } A_t^i\right)
\]
is a bounded martingale that goes to 0 when $t\rightarrow +\infty$. 
It follows that this martingale is identically zero almost surely, which implies $\phi=0$.
\end{enumerate}
\end{proof}
For later reference, we rephrase the second part of the previous proposition as follows: 
\begin{cor}\label{k}
The function $h_\mu^\theta(x)=e^{\mu(x)} j_\mu^\theta (x)$ is the unique solution to
\begin{equation}\label{pde2}
\frac{1}{2} \Delta h_\mu^\theta(x) - \sum_{i=1}^d \theta_i^2 e^{-2\alpha_i (x)} h_\mu^\theta(x)
= \frac{1}{2} \|\mu\|^2 h_\mu^\theta(x),
\end{equation}
such that $e^{-\mu(x)} h_\mu^\theta(x)$ is bounded and
\[\lim_{x\rightarrow \infty,x \in \Omega } e^{-\mu(x)} h_\mu^\theta(x)=1.\]
\end{cor}
\begin{example}
The following example has been widely studied (see, for example,~\cite{d,my} and references therein).
Suppose $n=d=1$, $\theta_1^2=1/2$ and $\alpha_1(x)=x$.  Then
\[A_\infty=\int_0^\infty e^{-2(B_t+\mu t)}dt , \quad \mu > 0\]
where $(B_t,\ t\ge 0)$ is a standard one-dimensional Brownian motion,
and \[j_\mu^\theta(x) = \E \left(\exp\left(-\frac{1}{2} e^{-2x} A_\infty\right)\right)\qquad x\in\R.\]
In this case, $h_\mu^\theta (x)=e^{\mu x}j_\mu^\theta(x) $ solves the equation
\[
\left( \frac{d^2}{dx^2} -e^{-2x} \right)h_\mu^\theta =\mu^2 h_\mu^\theta 
\]
This equation is  easily solved by means of Bessel functions.  By taking into account the boundary 
condition when $x \rightarrow +\infty$, we recover the formula~\cite[Theorem 6.2]{my}:
\begin{equation}\label{cf1}
j_\mu^\theta(x) = \frac{2^{1-\mu}}{\Gamma(\mu)} e^{-\mu x} K_\mu(e^{-x}),
\end{equation}
where $K_\mu$ is the Macdonald function~\cite{lebedev}:
\begin{equation}\label{mac}
K_\mu (x)=\frac{1}{2} \left(\frac{x}{2} \right)^\mu \int_0^{+\infty} \frac{e^{-t-\frac{x^2}{4t}} }{t^{1+\mu}} dt.
\end{equation}
The formula (\ref{cf1}) can also be derived using the fact~\cite{d} that 
$A_\infty$ has the same law as $1/2\gamma_\mu$, where $\gamma_\mu$ is
a gamma distributed random variable with parameter $\mu$.
\end{example}

\begin{example}
The following example has also been studied in the literature~\cite{grosche,im}.
Suppose $n=1$, $d=2$,  $\theta_1^2=\theta_2^2=1/2$,  $\alpha_1(x)=x$ and $\alpha_2 (x)=\frac{x}{2}$. 
Then
\[
A^1_\infty=\int_0^\infty e^{-2(B_t+\mu t)}dt, \quad A^2_\infty=\int_0^\infty e^{-(B_t+\mu t)}dt, \quad \mu > 0
\]
where $(B_t,\ t\ge 0)$ is a standard one-dimensional Brownian motion,
and \[j_\mu^\theta(x) = \E \left(\exp\left(-\frac{1}{2} e^{-2x} A^1_\infty -\frac{1}{2} e^{-x} A^2_\infty\right)\right)\qquad x\in\R.\]
In this case, $h_\mu^\theta (x)=e^{\mu x}j_\mu^\theta(x) $ solves the equation
\[
\left( \frac{d^2}{dx^2} -e^{-x}-e^{-2x} \right)h_\mu^\theta =\mu^2 h_\mu^\theta 
\]
This is Schr\"odinger's equation with the so-called Morse potential. 
It is solved by means of Whittaker functions and by taking into account the 
boundary condition when $x \rightarrow +\infty$, we get
\[
j_\mu^\theta(x)=2^{\mu-\frac{1}{2} }\frac{\Gamma(1+\mu)}{\Gamma (2\mu)} e^{\left( -\mu+\frac{1}{2} \right)x} W_{-\frac{1}{2},\mu} (2e^{-x})
\]
where $ W_{k,\mu}$ is the Whittaker function (see \cite[pp.279]{lebedev}):
\[
 W_{k,\mu} (x)=\frac{x^ke^{-\frac{x}{2}}}{\Gamma \left(\frac{1}{2}+\mu-k \right)}\int_0^{+\infty} e^{-t}t^{\mu-k-\frac{1}{2}} \left( 1+\frac{t}{x}\right)^{\mu+k-\frac{1}{2}} dt.
\]
\end{example}

\subsection{Conditional densities}

We prove now that the random variable $A_\infty$ has a smooth density with respect to the Lebesgue measure of $\mathbb{R}^d$ and morever give an expression of the conditional densities only in terms of this density.

\begin{proposition}\label{conditional}
The random variable $A_\infty$ has a smooth density $p$ with respect to the Lebesgue measure of $\mathbb{R}^d$ and for $t \ge 0$
\[
\mathbb{P} (A_\infty \in dy \left| \mathcal{F}_t \right.)=e^{2\sum_{i=1}^d \alpha_i (B^{(\mu)}_t)} 
p  \left( e^{2\alpha_1 (B^{(\mu)}_t)} (y_1-A^1_t),\cdots, e^{2\alpha_d (B^{(\mu)}_t)} (y_d-A^d_t)\right)  \mathbf{1}_{(0, y_1) \times \cdots \times (0, y_n)} (A_t)dy
\]
where $\mathcal{F}$ is the natural filtration of $B^{(\mu)}$. 
\end{proposition}
\begin{proof}
If we denote by $\phi$ the characteristic function of $A_\infty$:
\[
\phi (\lambda)=\mathbb{E} \left( e^{- (\lambda, A_\infty) }\right), \quad \lambda_1,...,\lambda_d >0,
\]
then, 
\begin{align*}
\mathbb{E} \left( e^{-  \sum_{i=1}^d \lambda_i A^i_\infty} \left| \mathcal{F}_t \right. \right) & =  e^{- \sum_{i=1}^d \lambda_i A^i_t}\mathbb{E} \left( e^{-  \sum_{i=1}^d \lambda_i (A^i_\infty-A^i_t)} \left| \mathcal{F}_t \right. \right)  \\
 & = e^{-(\lambda, A_t)} \phi \left( e^{-2\alpha_1 (B^{(\mu)}_t)} \lambda_1,\cdots,e^{-2\alpha_d (B^{(\mu)}_t)} \lambda_d \right).
\end{align*}
Therefore, the process $e^{-(\lambda ,A_t)} \phi \left( e^{-2\alpha_1 (B^{(\mu)}_t)} \lambda_1,\cdots,e^{-2\alpha_d (B^{(\mu)}_t)} \lambda_d \right)$ is a martingale. This implies that the function $e^{-(\lambda, a)} \phi \left( e^{-2\alpha_1 (x)} \lambda_1,\cdots,e^{-2\alpha_d (x)} \lambda_d \right)$ is harmonic for the operator $\frac{1}{2} \Delta_x +(\mu,\nabla_x)+\sum_{i=1}^d e^{-2\alpha_i (x)} \frac{\partial }{\partial a_i }$. This operator being hypoelliptic, this implies that $A_\infty$ has a smooth density with respect to the Lebesgue measure of $\mathbb{R}^d$. The result about the conditional densities stems from the injectivity of the Laplace transform.
\end{proof}

In particular, we deduce from the previous proposition that if for $y \in \mathbb{R}^d_+$, we denote
\[
q(x,a,y) =e^{2\sum_{i=1}^d \alpha_i (x)}p \left( e^{2\alpha_1 (x)}(y_1-a_1), \cdots, e^{2\alpha_d (x)}(y_d-a_d)\right), 
\]
for $0<a_i<y_i, x \in \mathbb{R}^d$,
then the process $q\left( B_t^{(\mu)}, A_t ,y \right)\mathbf{1}_{(0, y_1) \times \cdots \times (0, y_n)} (A_t)$ is a martingale. It implies that for any $y \in \mathbb{R}^d_+$, $q(x,a,y)$ satisfies the following partial differential equation:
\[
\frac{1}{2} \Delta_x q+(\mu,\nabla_x q)+\sum_{i=1}^d e^{-2\alpha_i (x)} \frac{\partial q}{\partial a_i }=0.
\]
It also implies that $p$ is a solution of the partial differential equation:
\begin{align*}
 \sum_{i,j=1}^d (\alpha_i,\alpha_j) y_i y_j \frac{\partial^2 p}{\partial y_i \partial y_j}
+\sum_{i=1}^d & \left( \left( \alpha_i (\mu) + \|\alpha_i \|^2+2\sum_{j=1}^d (\alpha_i,\alpha_j)\right)y_i-\frac{1}{2} \right)\frac{\partial p}{\partial y_i} \\
=& -\left( \sum_{i,j=1}^d (\alpha_i,\alpha_j)+\sum_{i=1}^d  \alpha_i (\mu)  \right) p.
\end{align*}

\begin{example}
Suppose $n=d=1$, $\theta_1^2=1/2$ and $\alpha_1(x)=x$.  Then $A_\infty$ is distributed as 
$1/2\gamma_\mu$, where $\gamma_\mu$ is a gamma law with parameter $\mu$, that is
\[
p(y)=\frac{1}{2^\mu \Gamma (\mu)} \frac{e^{-\frac{1}{2y}}}{y^{1+\mu}} 1_{\mathbb{R}_{>0}} (y),
\]
and we have
\[
q(x,a,y)= \frac{1}{2^\mu \Gamma (\mu)} \frac{e^{-2\mu x -\frac{1}{2}\frac{e^{-2x}}{y-a} } }{(y-a)^{1+\mu} }1_{\mathbb{R}_{>0}} (y-a)
\]
\end{example}

\begin{example}
Suppose $n=1$, $d=2$,  $\alpha_1(x)=x$ and $\alpha_2 (x)=\frac{x}{2}$. Then, as seen before,
\[
A^1_\infty=\int_0^\infty e^{-2(B_t+\mu t)}dt, \quad A^2_\infty=\int_0^\infty e^{-(B_t+\mu t)}dt, \quad \mu > 0,
\]
and for $\lambda_1,\lambda_2>0$,
\begin{align*}
  \mathbb{E} \left( e^{-\frac{1}{2} \lambda_1^2 A_1^\infty -\frac{1}{2} \lambda_2^2 A_2^\infty}\right) &
 = 2^{\mu- \frac{1}{2}} \lambda_1^{\mu- \frac{1}{2}} \frac{ \Gamma\left(\mu+\frac{1}{2}+\frac{\lambda_2^2}{2\lambda_1} \right)}{\Gamma(2\mu)}  W_{-\frac{\lambda_2^2}{2\lambda_1},\mu} (2\lambda_1) \\
 &= \frac{e^{-\lambda_1}}{\Gamma (2 \mu)} \int_0^{+\infty} e^{-t} t^{\mu+\frac{\lambda_2^2}{2\lambda_1}-\frac{1}{2}} (2\lambda_1+t)^{\mu-\frac{\lambda_2^2}{2\lambda_1}-\frac{1}{2}} dt.
\end{align*}
By using in the previous integral the change of variable $t=\frac{2\lambda_1}{e^{\lambda_1 s}-1}$, we deduce the following nice formula
\[
\mathbb{E} \left( e^{ -\frac{1}{2} \lambda_1^2 A^1_\infty} \mid A_\infty^2 =s \right) \mathbb{P} (A^2_\infty \in ds)
=\frac{\lambda_1^{2\mu+1}}{2\Gamma(2\mu)} \frac{e^{-\lambda_1 \mathrm{cotanh} \frac{\lambda_1 s}{2}}}{\left( \sinh \frac{\lambda_1 s}{2}\right)^{2\mu+1} }ds, \quad s>0.
\]
This conditional Laplace transform can be inverted (see for instance \cite{ghomrasni}) but, unlike the one-dimensional case, it does not seem to lead to a nice formula for $p$:
\begin{align*}
p(y_1,y_2)=\frac{2^{2\mu}}{\Gamma (2\mu) \sqrt{2\pi}} \sum_{j,k=0}^{+\infty}&  \frac{(-1)^j 2^j}{j!} \frac{\Gamma(j+2\mu+1+k)}{k!\Gamma(j+2\mu+1)} 
\frac{1}{y_1^{\frac{j}{2}+\mu +\frac{3}{2}}} e^{-\frac{(1+y_2(k+j+\mu+\frac{1}{2}))^2}{4y_1}}  \\
& \times D_{j+2\mu+2} \left(\frac{1+y_2(k+j+\mu+\frac{1}{2})}{\sqrt{y_1}} \right)
\end{align*}
where $D_\nu$ is the parabolic cylinder function such that
\[
\int_0^{+\infty} \frac{e^{-\theta t}}{t^{1+\nu}} e^{-\frac{a^2}{4t}} D_{2\nu+1} \left( \frac{a}{\sqrt{t}}\right)dt
=\sqrt{\pi} 2^{\nu+\frac{1}{2}} \theta^{\nu} e^{-a\sqrt{2\theta}},
\]
that is
\[
D_\nu (x) =\frac{\sqrt{2}}{\sqrt{\pi}} e^{\frac{x^2}{4}}
\int_0^{+\infty} t^\nu e^{-\frac{t^2}{2}} \cos \left( xt-\frac{\pi
\nu}{2} \right) dt, \quad \nu >-1
\]
\end{example}

\section{Brownian motion conditioned on its exponential functionals}

In this section, we study the Doob transforms of the process $\left( B_t^{(\mu)} ,A_t \right)$ 
associated with the conditioning of $A_\infty$. 
We first start with the bridges which are the extremal points.

\begin{lem}[Equation of the bridges]
Let $y \in \mathbb{R}_+^d$. The law of the process $(B_t+\mu t)_{t \ge 0}$ conditioned by
\[
A_{\infty} =y
\]
solves  the following stochastic differential equation:
\[
dX_t=\left(\mu+(\nabla_x \ln q)\left(X_t , \int_0^te^{-2\alpha(X_s)}ds,y\right)\right) dt+d\beta_t
\]
where, $(\beta_t)_{t\ge 0}$ is a standard Brownian motion 
\end{lem}

\begin{proof}
This follows directly from Proposition \ref{conditional} and Girsanov's theorem.
\end{proof}

\begin{example}
The following example is considered~\cite{my3}.
Suppose $n=d=1$, $\theta_1^2=1/2$ and $\alpha_1(x)=x$. 
Then the equation becomes
\[
dX_t=\left(-\mu+\frac{e^{-2X_t}}{y-\int_0^t e^{-2X_s} ds}\right) dt+d\beta_t .
\]
\end{example}

Let $\mathbb{P}^\mu$ be the law of $B^{(\mu)}$ and $\pi$ be the coordinate process on the space of continuous functions $\mathbb{R}_+ \rightarrow \mathbb{R}^n$.
If $\nu$ is a probability measure on $\mathbb{R}^d_+$, in what follows (see \cite{b2}), we call the probability 
\[
\int_{\mathbb{R}^d_+} \mathbb{P}^\mu \left( \quad \cdot \quad \left| \int_0^{+\infty} e^{-2\alpha_1 (\pi_s)}ds=y_1,\cdots, \int_0^{+\infty} e^{-2\alpha_d (\pi_s) }ds=y_d \right) \right. \nu (dy),
\]
the law of the process $(B_t+\mu t)_{t \ge 0}$ conditioned by
\[
A_{\infty} =^{law} \nu .
\]

\begin{prop}
Let $v$ be a bounded and positive function such that $\int_{\mathbb{R}^d} v(y)p(y)dy =1$. The law of the process $(B_t+\mu t)_{t \ge 0}$ conditioned by
\[
A_{\infty} =^{law} v(x) p(x) dx
\]
solves  the following stochastic differential equation:
\[
dX_t=\left(\mu+F_v \left( \int_0^t e^{-2\alpha_1 (X_s)} ds, \cdots, \int_0^t e^{-2\alpha_d (X_s)} ds, X_t \right)\right)dt+d\beta_t
\]
where, $(\beta_t)_{t\ge 0}$ is a standard Brownian motion and $F_v: \mathbb{R}^d \times \mathbb{R}^n \rightarrow \mathbb{R}^n$ is given by
\[
F_v(a,x)=(\nabla_x \ln \phi_v) (a,x)
\]
with
\[
\phi_v (a,x)=\int_{\mathbb{R}^d} p(z) v\left(a_1+e^{-2\alpha_1 (x)}z_1,\cdots,a_d+e^{-2\alpha_d (x)}z_d  \right) dz.
\]
\end{prop}

\begin{proof}

Following \cite{b2},  we have to write the stochastic differential equation associated with the conditioning
\[
A_{\infty} =^{law} p(x) v(x) dx.
\]
But
\begin{align*}
  \mathbb{E} \left( v\left( A_{\infty} \right)  \left| \mathcal{F}_t \right. \right) 
=&e^{2\sum_{i=1}^d \alpha_i (B^{(\mu)}_t)} \int_{\mathbb{R}^d} p\left( e^{2\alpha_1 (B^{(\mu)}_t)} (y_1-A^1_t),\cdots, e^{2\alpha_d (B^{(\mu)}_t)} (y_d-A^d_t)\right)v(y)dy \\
=&\phi_v (A_t,B_t^{(\mu)}),
\end{align*}
so that we get the expected conditioned stochastic differential equation by Girsanov theorem.

\end{proof}

In the previous proposition, the drift $F_v (a,x)$ depends only on $x$ if, and only if, 
\[v(x)=\frac{e^{-\sum_{i=1}^d \theta_i^2 x_i}}{j_\mu^\theta (0)},\] 
for some $\theta \in \mathbb{R}^d$.
Therefore: 
\begin{cor}\label{loi}
For $\theta \in \mathbb{R}^d$, the law of the process $(B_t+\mu t)_{t \ge 0}$ conditioned by
\[
A_{\infty} =^{law} \frac{e^{-\sum_{i=1}^d \theta_i^2 x_i}}{j_\mu^\theta (0)} p(x) dx
\]
is the law of a Markov process.
Moreover, in that case, it solves in law the following stochastic differential equation
\begin{equation}\label{conditioning}
dX_t=\nabla \ln h_\mu^\theta (X_t)dt+d\beta_t.
\end{equation}
\end{cor}

We now show that the  pathwise uniqueness property holds for the stochastic differential equation (\ref{conditioning}). In what follows, we denote
\[
\mathcal{L}_\mu^{\theta}=\nabla \ln h_\mu^\theta \nabla+\frac{1}{2} \Delta.
\]
Let us observe that for the generator $\mathcal{L}_\mu^{\theta}$, we have a useful intertwining with the Schr\"odinger operator $\frac{1}{2} \Delta  - \sum_{i=1}^d \theta_i^2 e^{-2\alpha_i (x)} - \frac{1}{2} \|\mu\|^2$ that will be used several times in the sequel.
\begin{proposition}
\[
h_\mu^\theta \mathcal{L}_\mu^{\theta} =\left( \frac{1}{2} \Delta  - \sum_{i=1}^d \theta_i^2 e^{-2\alpha_i (x)} - \frac{1}{2} \|\mu\|^2\right) h_\mu^\theta
\]
\end{proposition}

\begin{proof}
If $f$ is a smooth function then we have
\begin{align*}
 & \left( \frac{1}{2} \Delta  - \sum_{i=1}^d \theta_i^2 e^{-2\alpha_i (x)} - \frac{1}{2} \|\mu\|^2\right) (h_\mu^\theta f) \\
 =& \frac{1}{2} (\Delta h^\theta_\mu) f +\frac{1}{2} (\Delta f) h^\theta_\mu+ \nabla f \nabla h^\theta_\mu  - \left( \sum_{i=1}^d \theta_i^2 e^{-2\alpha_i (x)} +\frac{1}{2} \|\mu\|^2 \right)(h_\mu^\theta f).
\end{align*}
Since
\[
\frac{1}{2} \Delta h^\theta_\mu=\left(\sum_{i=1}^d \theta_i^2 e^{-2\alpha_i (x)} +\frac{1}{2} \|\mu\|^2 \right) h_\mu^\theta,
\]
the result readily follows.
\end{proof}
We can now deduce:
\begin{theorem}
Let $\theta \in \mathbb{R}^d$. If $(\beta_t)_{t \ge 0}$ is a Brownian motion, then for $x_0 \in \mathbb{R}^n$, there exists a unique process
 $(X_t^{x_0})_{t\ge0}$ adapted to the filtration of $(\beta_t)_{t \ge 0}$ such that:
 \begin{align}\label{EDS}
 X_t^{x_0}=x_0+\int_0^t \nabla \ln h_\mu^\theta (X^{x_0}_s)ds+\beta_t, \quad t \ge 0.
 \end{align}
Moreover, in law, the process $(X_t^{x_0})_{t\ge0}$ is equal to 
$(B_t +\mu t +x_0)_{t \ge 0}$ conditioned by:
\begin{align*}
& \left(  \int_0^{+\infty} e^{-2\alpha_i (B_t +\mu t+x_0) }dt \right)_{1 \le i \le d} \\
& =^{law}\frac{e^{- \sum_{i=1}^d \theta_i^2  y_i}e^{2\sum_{i=1}^d \alpha_i (x_0) } 
p ( e^{2\alpha_1 (x_0) } y_1, \cdots,e^{2\alpha_d (x_0) } y_d ) }{j_\mu^\theta (x_0) } dy .
\end{align*}
\end{theorem}

\begin{proof}
Let $x_0 \in \mathbb{R}^n$. Since the function $\nabla \ln h_\mu^\theta$ is locally Lipschitz, up to an 
explosion time $\mathbf{e}$ we have a unique solution $ X_t^{x_0}$ for the equation (\ref{EDS}). Our goal is now to show
that almost surely $\mathbf{e}=+\infty$. For that, we construct a suitable Lyapunov function for the generator
$\mathcal{L}_\mu^{\theta}$.

Let
\[
U(x)=\frac{\cosh 2(\mu,x)}{h^\theta_\mu (x)}.
\]
It is easily seen that when $\parallel x \parallel \rightarrow +\infty$, $U(x) \rightarrow +\infty$. Moreover, from the
intertwining,
\[
h^\theta_\mu \mathcal{L}_\mu^{\theta} U=
\left( \frac{1}{2} \Delta  - \sum_{i=1}^d \theta_i^2 e^{-2\alpha_i (x)} - \frac{1}{2} \|\mu\|^2\right) \cosh 2(\mu ,x)
\le \frac{3}{2} \| \mu \|^2 \cosh 2(\mu ,x).
\]
Therefore
\[
\mathcal{L}_\mu^{\theta} U \le  \frac{3}{2} \| \mu \|^2 U
\]
It implies that the process $(e^{- \frac{3}{2} \| \mu \|^2  t \wedge \mathbf{e}} U (X^{x_0}_{t \wedge \mathbf{e}}))_{t \ge 0}$ is a positive
supermartingale. Since $U(x) \rightarrow +\infty$ when $\parallel x \parallel \rightarrow +\infty$, we deduce that
almost surely $\mathbf{e}=+\infty$.

Consequently, there is a unique solution $ (X_t^{x_0})_{t \ge 0}$ for the equation (\ref{EDS}). The second part of the 
theorem is a direct consequence of Corollary \ref{loi} and  uniqueness in law for the equation (\ref{EDS}).  
\end{proof}

\begin{example}
Suppose $n=d=1$ and $\theta_1^2=1/2$ and $\alpha_1(x)=x$. Then 
\begin{equation}\label{myd}
\mathcal{L}_\mu^{\theta}=\left(\mu+e^{-x}\frac{K_{\mu-1} (e^{-x})}{K_\mu (e^{-x})} \right)\frac{d}{dx}+\frac{1}{2} \frac{d^2}{dx^2}
\end{equation}
Let us denote $p^{\mu,\theta}_t (x,y)$ the heat kernel of $\mathcal{L}_\mu^{\theta}$. From the intertwining, we have
\[
p^{\mu,\theta}_t (x,y)=\frac{K_\mu (e^{-y})}{K_\mu (e^{-x})} q^{\mu,\theta}_t (x,y),
\]
where $q^{\mu,\theta}_t (x,y)$ is the heat kernel of $\frac{1}{2} \left( \frac{d^2}{dx^2} -e^{-2x}-\mu^2\right)$. 
This kernel can be explicitly computed (see \cite{ams} or \cite[Remark 4.1]{my2}):
\[
q^{\mu,\theta}_t (x,y)=e^{-\frac{\mu^2}{2}t} \int_0^{+\infty} \exp \left(-\frac{\xi}{2} -\frac{e^{-2x}+e^{-2y}}{2\xi} \right)
\Theta \left( \frac{ e^{-x-y}}{\xi}, t \right) \frac{d\xi}{\xi},
\]
with
\[
\Theta (r,t)=\frac{r}{\sqrt{2\pi^3 t}} e^{\frac{\pi^2}{2t}} \int_0^{+\infty} e^{-\frac{\xi^2}{2t}} e^{-r \cosh \xi} \sinh \xi \sin \frac{\pi \xi}{t} d\xi.
\]
We deduce from that
\[
p^{\mu,\theta}_t (-\infty,y)=2 e^{-\frac{\mu^2 t}{2}} \Theta (e^{-y},t)K_\mu (e^{-y}),
\]
so that $-\infty$ is an entrance point for the diffusion with generator $\mathcal{L}_\mu^{\theta}$.

The resolvent kernel of $(-\mathcal{L}_\mu^{\theta}+\frac{\alpha^2}{2})^{-1}$ 
is also easily computed:
\[G^{\mu,\theta} (x,y,-\frac{\alpha^2}{2})=2 \frac{K_\mu (e^{-y})}{K_\mu (e^{-x})} I_{\sqrt{\alpha^2+\mu^2}} (e^{-y}) K_{\sqrt{\alpha^2+\mu^2}} (e^{-x}), \quad x \le y. 
\]
And we can observe that
\[
G^{\mu,\theta} (-\infty,y,-\frac{\alpha^2}{2})=2K_\mu (e^{-y})I_{\sqrt{\alpha^2+\mu^2}} (e^{-y}).
\]
\end{example}

\begin{example}
Suppose $n=1$, $d=2$,  $\theta_1^2=\theta_2^2=1/2$,  $\alpha_1(x)=x$ and $\alpha_2 (x)=\frac{x}{2}$. In that case
\[
\mathcal{L}_\mu^{\theta}=\left(\frac{1}{2}-2e^{-x} \frac{W'_{-1/2,\mu}(2e^{-x})}{W_{-1/2,\mu}(2e^{-x})}  \right)\frac{d}{dx}+\frac{1}{2} \frac{d^2}{dx^2}
\]
and
\[
p^{\mu,\theta}_t (x,y)=e^{\frac{1}{2} (y-x) }\frac{W_{-\frac{1}{2},\mu} (2e^{-y})}{W_{-\frac{1}{2},\mu} (2e^{-x})} q^{\mu,\theta}_t (x,y),
\]
where $q^{\mu,\theta}_t (x,y)$ is the heat kernel of $\frac{1}{2} \left( \frac{d^2}{dx^2} -e^{-x}-e^{-2x}-\mu^2\right)$. 
We have (see \cite{ams} or \cite[p342]{my2}):
\[
q^{\mu,\theta}_t (x,y)=e^{-\frac{\mu^2}{2}t} \int_0^{+\infty} e^{-\xi -(e^{-x}+e^{-y}) \mathrm{cotanh} \xi } \Theta \left(2 \frac{e^{-\frac{x+y}{2}}}{\sinh \xi}, \frac{t}{\xi} \right) \frac{d\xi}{\sinh \xi}.
\]

The resolvent kernel of $(-\mathcal{L}_\mu^{\theta}+\frac{\alpha^2}{2})^{-1}$, is for $x \le y$:
\[
G^{\mu,\theta} (x,y,-\frac{\alpha^2}{2})=\frac{\Gamma \left(1+ \sqrt{\alpha^2+\mu^2}\right) }{ \Gamma \left(1+ 2\sqrt{\alpha^2+\mu^2} \right) } \frac{W_{-\frac{1}{2},\mu} (2e^{-y})}{W_{-\frac{1}{2},\mu} (2e^{-x})} 
W_{-\frac{1}{2},\sqrt{\alpha^2+\mu^2}} (2e^{-x})M_{-\frac{1}{2},\sqrt{\alpha^2+\mu^2}} (2e^{-y}), 
\]
and we get:
\[
G^{\mu,\theta} (-\infty,y,-\frac{\alpha^2}{2})=\frac{\Gamma \left(1+ \sqrt{\alpha^2+\mu^2}\right) }{ \Gamma \left(1+ 2\sqrt{\alpha^2+\mu^2} \right) }
W_{-\frac{1}{2},\mu} (2e^{-y}) M_{-\frac{1}{2},\sqrt{\alpha^2+\mu^2}} (2e^{-y}),
\]
so that $-\infty$ is also an entrance point for the diffusion with generator $\mathcal{L}_\mu^{\theta}$.
\end{example}

Motivated by the two previous examples, the question of existence of entrance laws for the diffusion
with generator $\mathcal{L}_\mu^{\theta}$ is natural. As a general result, we can prove:
\begin{proposition}
Assume $n=1$, $\alpha_1,\cdots,\alpha_d >0$ and $\theta \in \mathbb{R}^d-\{0\}$, then $-\infty$ is an entrance point for the diffusion with generator  $\mathcal{L}_\mu^{\theta}$.
\end{proposition}
\begin{proof}
Without loss of generality, we can assume that $\theta_1 >0$.
Let us recall $h_\mu^\theta$ solves the Schr\"odinger equation
\[
\frac{1}{2} (h_\mu^\theta)'' =\left( \sum_{i=1}^d \theta_i^2 e^{-2\alpha_i x} +\frac{1}{2} \mu ^2 \right)h_\mu^\theta,
\]
and that $k_\mu^{\theta_1}(x)=K_\mu (\theta_1 e^{-\alpha_1x})$ solves the equation:
\[
\frac{1}{2} (k_\mu^{\theta_1})'' =\left( \theta_1^2 e^{-2\alpha_1 x} +\frac{1}{2} \mu ^2 \right)k_\mu^{\theta_1}.
\]
Let $W(x)=k_\mu^{\theta_1}(x)(h_\mu^{\theta})'(x)-(k_\mu^{\theta_1})'(x)(h_\mu^{\theta})(x)$. Since
\[
W'(x)=k_\mu^{\theta_1}(x)(h_\mu^{\theta})''(x)-(k_\mu^{\theta_1})''(x)(h_\mu^{\theta})(x)\ge 0,
\]
we deduce that $W$ is increasing. Moreover, it is easily seen that $\lim_{x \to -\infty} W(x)=0$. Therefore $W \ge 0$.
Hence $\frac{(h_\mu^\theta)'}{h_\mu^\theta} (x) \ge 
-\alpha_1 \theta_1 e^{-\alpha_1 x} \frac{K'_\mu (\theta_1 e^{-\alpha_1x})}{K_\mu (\theta_1 e^{-\alpha_1x})}$.

Now, from the comparison principle for stochastic differential equations, we deduce that if, for $x \in \mathbb{R}$, we denote
$(X_t^x)_{t \ge 0}$ and $(Y_t^x)_{t \ge 0}$ the solutions  of the stochastic differential equations,
\[ 
X_t^x=x+\int_0^t \frac{(h_\mu^\theta)'}{h_\mu^\theta}(X_s^x) ds +\beta_t
\]
\[
Y_t^x=x+\int_0^t -\alpha_1 \theta_1 e^{-\alpha_1 Y_s^x} 
\frac{K'_\mu (\theta_1 e^{-\alpha_1 Y_s^x})}{K_\mu (\theta_1 e^{-\alpha_1 Y_s^x})} ds +\beta_t
\]
where $(\beta_t)_{t \ge 0}$ is a standard Brownian motion, then we have almost surely
\[
X_t^x \ge Y_t^x.
\]
Since $-\infty$ is an entrance point for the diffusion $(Y_t^x)_{t \ge 0, x \in \mathbb{R}}$,
we deduce that  $-\infty$ is an entrance point for the diffusion with generator  $\mathcal{L}_\mu^{\theta}$ .

\end{proof}
We conjecture the existence of entrance laws for $n \ge 1$, but let us observe that, in general, we do not have unicity.
Indeed, let us consider the following example
\[
n=2, d=1, \alpha (x)=\frac{x_2-x_1}{\sqrt{2}}, \theta^2=\frac{1}{2}.
\]
In that case, by using one dimensional results, we compute:
\[
h_\mu^\theta (x)= \frac{2^{1-\alpha(\mu)}}{\Gamma (\alpha (\mu))} e^{\alpha^*(\mu)\alpha^*(x)}K_{\alpha (\mu)} (e^{-\alpha(x)}),
\]
where $\alpha^*(x)=\frac{x_2+x_1}{\sqrt{2}}$.
The heat kernel of $\mathcal{L}_\mu^{\theta}$ is also explicitly given by
\begin{eqnarray*}
\lefteqn{
p^{\mu,\theta}_t (x,y)=e^{-\frac{1}{2} \| \mu \|^2 t} \frac{h_\mu^\theta (y)}{h_\mu^\theta (x)}
\frac{1}{\sqrt{2\pi t}} e^{-\frac{(\alpha^* (x)-\alpha^*(y))^2}{2t}} } \\
&& \times \int_0^{+\infty} \exp \left(-\frac{\xi}{2} 
-\frac{e^{-2\alpha(x)}+e^{-2\alpha(y)}}{2\xi} \right)
\Theta \left( \frac{ e^{-\alpha(x)-\alpha(y)}}{\xi}, t \right) \frac{d\xi}{\xi}.
\end{eqnarray*}
And we deduce that when $\alpha(x) \to -\infty$ with $\alpha^* (x) \to k \in \mathbb{R}$,
\[
p^{\mu,\theta}_t (x,y) \to 2 e^{-\frac{1}{2} \| \mu \|^2 t} h_\mu^\theta (y) e^{-k \alpha^* (\mu)} \frac{1}{\sqrt{2\pi t}} 
e^{-\frac{(k-\alpha^*(y))^2}{2t}} \Theta (e^{-\alpha(y)},t).
\]
Therefore, in that case we get an infinite set of entrance laws when $\alpha (x) \to - \infty$.

\section{Whittaker functions}

From now on we consider the case where $\Pi=\{\a_1,\ldots,\a_d\}$ is a
simple system.  In other words:
\begin{enumerate}
\item  The vectors $\a_1,\ldots,\a_d$ are linearly independent;
\item the group $W$ generated by reflections through the hyperplanes 
\[H_\a=\{x\in\R^n:\ \a(x)=0\},\qquad \a\in\Pi\] is finite;
\item $\{x\in\R^n:\ \a(x)\ge 0,\ \forall\a\in\Pi\} $ is a
fundamental domain for the action of $W$ on $\R^n$;
\item $2(\a,\beta)/(\a,\a)\in\Z$ for all $\a,\beta\in\Pi$.
\end{enumerate}
In this setting, the Schr\"odinger operator
\[
H=\frac{1}{2} \Delta  - \sum_{i=1}^d \theta_i^2 e^{-2\alpha_i (x)} 
\]
is the Hamiltonian of the (generalized) quantum Toda lattice (see, for example, \cite{an}). 
The function $h_\mu^\theta$ considered in the previous section can be expressed in terms 
of a particular eigenfunction of $H$, known as a class-one Whittaker function.

\subsection{Class one Whittaker functions}
Class one Whittaker functions associated with semisimple Lie groups were 
introduced by Kostant~\cite{k} and Jacquet~\cite{j}, and have been studied extensively in the literature. 
They are closely related to Whittaker models of principal series representations
and play an important role in the study of automorphic forms associated with Lie 
groups~\cite{bump}. 
They also arise as eigenfunctions of the (generalised) quantum 
Toda lattice~\cite{k,an}.  
For completeness we will describe briefly the abstract definition of class-one Whittaker functions, 
following~\cite{h}.  

Let $G$ be a connected, noncompact, semisimple Lie group with finite centre.
Let $\mathfrak g_0$ be the Lie algebra of $G$ with complexification $\mathfrak g$.
Denote by $B(\cdot,\cdot)$ the Killing form on $\mathfrak g$.
Let $K$ be a maximal compact subgroup of $G$ with Lie algebra $\mathfrak k_0$
and denote the complexification of $\mathfrak k_0$ by $\mathfrak k$.
Let $\mathfrak p_0$ be the orthogonal complement of $\mathfrak k_0$ in
$\mathfrak g_0$ with respect to the Killing form.  Let $\theta$ be the corresponding
Cartan involution.  Let $\mathfrak a_0$ be a maximal abelian subspace in $\mathfrak p_0$
and denote its complexification by $\mathfrak a$.
Denote by $\Sigma$ the set of all non-zero roots of $\mathfrak g_0$ relative to $\mathfrak a_0$.  
For $\alpha\in\Sigma$, denote by $m(\alpha)$ the dimension of the root space 
\[\mathfrak g_0^\alpha=
\{ X\in \mathfrak g_0:\ \mbox{ad} (H) X=\alpha(H) X\mbox{ for all } H\in \mathfrak a_0\}.\]
Let $\Sigma_+$ be a positive system of roots in $\Sigma$ and let $\Pi=\{\alpha_1,\ldots,\alpha_d\}$
be the corresponding set of simple roots.  Let 
$\mathfrak n_0= \sum_{\alpha\in\Sigma_+} \mathfrak g_0^\alpha$
and $N=\exp(\mathfrak n_0)$.  Then $G=NAK$ is an Iwasawa decomposition of $G$.
Let $\psi$ be a non-degenerate (unitary) character of $N$.  
Let $\eta$ be the unique Lie algebra homomorphism of $\mathfrak n_0$ into $\R$ such that
$\psi(n)=\exp(i\eta(X))$ for $n=\exp(X)\in N$.  
For each $\alpha\in\Sigma_+$,
let $X_{\alpha,i}\ (1\le i\le m(\alpha))$ be a basis of $\mathfrak g_0^\alpha$ satisfying 
$B(X_{\alpha,i},\theta X_{\alpha,j})=-\delta_{ij}\ (1\le i,j \le m(\alpha))$.
Denote by $\eta_\alpha$ the restriction of $\eta$ to $\mathfrak g_0^\alpha$
and set $|\eta_\alpha|^2=\sum_{1\le i\le m(\alpha)} \eta(X_{\alpha,i})^2$.
Denote by $U(\mathfrak g)$ and $U(\mathfrak a)$ the universal enveloping algebras of 
$\mathfrak g$ and $\mathfrak a$, respectively.  Let $\gamma$ denote the Harish-Chandra 
homomorphism from $U(\mathfrak g)^{\mathfrak k}$, the centraliser of $\mathfrak k$ in 
$U(\mathfrak g)$, into $U(\mathfrak a)$. For $\nu\in{\mathfrak a}^*$ and 
$z\in U(\mathfrak g)^{\mathfrak k}$, define $\chi_\nu(z)=\gamma(z)(\nu)$.  
The space of Whittaker functions on $G$ associated with $\nu\in{\mathfrak a}^*$, 
denoted $C^\infty_\psi(G/K,\chi_\nu)$, is the space of smooth functions on $G$ 
which satisfy:
\begin{enumerate}
\item $f(ngk)=\psi(n)f(g)$ for $n\in N$, $g\in G$ and $k\in K$, and 
\item $zf=\chi_\nu(z)f$ for $z\in U(\mathfrak g)^{\mathfrak k}$. 
\end{enumerate}
Set $\rho=\frac{1}{2} \sum_{\alpha\in\Sigma_+}m(\alpha)\alpha$.
For $g\in G$, define 
$1_\nu(g)=h(g)^{\nu+\rho}$
where $g=n(g)h(g)k(g)$ is the Iwasawa decomposition of $g$.  Let $s_0$ be the longest
element in $W$.  The {\em class-one Whittaker function}
associated with $\nu\in\mathfrak a^*$ is defined by
\begin{equation}\label{cow}
W_\nu(g)=\int_N 1_\nu(s_0 ng)\psi^{-1}(n) dn, \qquad g\in G.
\end{equation}
The convergence of this integral was established by Jacquet~\cite{j}.
For $\nu\in \mathfrak a^*$ and $\alpha\in\Sigma$, write $\nu_\a=(\a,\nu)/(\a,\a)$.
Let
\[D=\{\nu\in\mathfrak a^*:\ \Re(\nu_\alpha)>0,\mbox{ for all }\alpha\in\Sigma_+\}.\]
We record the following lemma for later reference.
\begin{lem}\label{bounded}
Let $\nu\in D$.  Then $h^{-s_0\nu-\rho} W_\nu(h)$ is uniformly bounded for $h\in A$.
\end{lem}
\begin{proof}
Gindikin and Karpelevich~\cite{gk} proved that the integral
\[c(\nu)=\int_N 1_\nu(s_0 n) dn\]
is absolutely convergent.
From (\ref{cow}) we can write 
\[W_\nu(h)=h^{s_0\nu+\rho} \int_N 1_\nu(s_0 n)\psi^{-1}(hnh^{-1}) dn \qquad h\in A.\]
Since $\psi$ is unitary, it follows that $h^{-s_0\nu-\rho} W_\nu(h)$ is bounded,
as required.
\end{proof}
\begin{rem}
In the above, $c(\nu)$ is the Harish-Chandra $c$-function.
\end{rem}

\subsection{Fundamental Whittaker functions}
Since $W_\nu(nhk)=\psi(n)W_\nu(h)$, all of the important information about
$W_\nu$ is contained in its restriction to $A$.  This leads to a more concrete description
which can be presented entirely in the context of the root system $\Sigma$.  
Readers not familiar with root systems may find it helpful to think of the `type $A$' case,
for example if $G=SL(n,\R)$. In this case, we can identify $\mathfrak a_0$ (and its dual) with 
\[\R_0^n=\{\l\in\R^n:\ \l_1+\cdots+\l_n=0\},\] 
and take $\Sigma=\{e_i - e_j,\ i\ne j\}$, 
$\Sigma_+=\{e_i-e_j,\ 1\le i<j\le n\}$ and $\Pi=\{e_i-e_{i+1},\ 2\le i\le n\}$, where
$\{e_1,\ldots,e_n\}$ is the standard basis for $\R^n$.
In general, the root system $\Sigma$ is 
{\em crystallographic}, that is, the numbers $2(\alpha,\beta)/(\alpha,\alpha),\ \alpha,\beta\in\Pi$ 
are all integers, and the $\Z$-span of $\Pi$ is a regular lattice in $\mathfrak a_0^*$. 
Since the Killing form is positive definite on $\mathfrak a_0^*$, it induces an inner product
$(\cdot,\cdot)$ on $\mathfrak a_0^*$, which extends to a nondegenerate bilinear form on 
$\mathfrak a^*$.
The following construction is due to Hashizume~\cite{h}.
Consider the lattice $L= 2\Z_+(\Pi)$, and set
$'\mathfrak a^*=\{\nu\in \mathfrak a^*:\ (\l,\l)+2(\l,\nu)\ne 0,\ \forall\lambda\in L\backslash\{0\} \}.$
For each $\nu\in {'\mathfrak a^*}$, define a set of real numbers $\{c_\l(\nu), \l\in L\}$ 
recursively as follows. Set $c_0(\nu)=1$ and 
\begin{equation}\label{recursion}
((\l,\l)+2(\l,\nu))c_\l(\nu)=2\sum_\a |\eta_\a|^2 c_{\l-2\a}(\nu) \qquad \l\in L,
\end{equation}
with the convention that $c_\lambda(\nu)=0$ if $\lambda\notin L$.
In~\cite{h} it is shown that the series 
\[\Phi_\nu(x) = \sum_{\l\in L} c_\l(\nu) e^{-(\l+\nu)(x)},\]
converges absolutely and uniformly for $x\in\mathfrak a$ and $\nu\in {'\mathfrak a^*}$.
Define $U$ to be the set of $\nu\in {'\mathfrak a^*}$ such that:
\begin{enumerate}
\item $\nu_\alpha\ne 0$ for all $\alpha\in\Sigma$;
\item $s\nu\in\ {'\mathfrak a^*}$ for all $s\in W$;
\item $s\nu-t\nu\notin\sum_{\alpha\in\Pi}\Z\alpha$ for any pair $s,t\in W$ such that $s\ne t$.
\end{enumerate}
For $s\in W$ denote by $l(s)$ the length of $s$. 
For $\nu\in U$, define $M(s,\nu)\ (s\in W)$, recursively as follows.   
For $s=s_\alpha\ (\a\in\Pi)$,
\[M(s_\a,\nu)=\left(|\eta_\a|/2\sqrt{2(\a,\a)}\right)^{2\nu_\a} e_\alpha(\nu) e_\alpha(-\nu)^{-1},\]
where
\[e_\alpha(\nu)^{-1}=\G((\nu_\a+m(\alpha)/2+1)/2)\G((\nu_\a+m(\alpha)/2+m(2\alpha))/2) .\]
If $s\in W$ and $\a\in\Pi$ such that $l(s_\a s)=l(s)+1$, then
\[M(s_\a s,\nu)=M(s,\nu)M(s_\a ,s\nu).\]
Let $\Sigma_+^\circ$ be the set of $\alpha\in\Sigma_+$ such that $\alpha/2$ is not a root.
The Harish-Chandra $c$-function is given by 
\[c(\nu)= \prod_{\a\in\Sigma_+^\circ} d_\alpha f_\a(\nu),\]
where \[f_\a(\nu)= \frac{\G(\nu_\a)\G((\nu_\a+m(\a)/2)/2)}
{\G(\nu_\a+m(\a)/2)\G((\nu_\a+m(\a)/2+m(2\a))/2)},\]
and
\[d_\a = 2^{(m(\a)-m(2\a))/2} (\pi/(\a,\a))^{(m(\a)-m(2\a))/2}.\]
Now define, for $\nu\in U$,
\begin{equation}\label{psi}
\Psi_\nu(x)=\sum_{s\in W} M(s_0s,\nu)c(s_0s\nu)\Phi_{s\nu}(x).
\end{equation}
Observe that $\Psi_\nu$ satisfies the functional equations
\begin{equation}\label{fe}
\Psi_\nu(x)=M(s,\nu)\Psi_{s\nu}(x)\qquad s\in W.
\end{equation}
Although the above construction places a restriction on $\nu$,
it is known that, for each $x\in\mathfrak a$, $\Psi_\nu(x)$ can be extended to
an entire function of $\nu\in \mathfrak a^*$.  In~\cite{h} it is shown that, for $x\in\mathfrak a$, 
$W_\nu(e^{-x})=e^{-\rho(x)}\Psi_\nu(x)$, so that 
\[W_\nu(g)= \psi(n(g)) h(g)^\rho\Psi_\nu(\log h(g)).\]
The functions $V_\nu$ defined by
\[V_\nu(g)= \psi(n(g)) h(g)^\rho\Phi_\nu(\log h(g)),\]
are called {\em fundamental} Whittaker functions.  
In~\cite{h} it is also shown that, for each $\nu\in U$, $\{V_{s\nu},\ s\in W\}$ form a 
basis for $C^\infty_\psi(G/K,\chi_\nu)$.

\subsection{The quantum Toda lattice}

As observed by Kostant~\cite{k}, Whittaker functions are eigenfunctions for the (generalised)
quantum Toda lattice. 
Denote by $\Delta$ the Laplacian on $\mathfrak a_0$ corresponding to the Killing form.
For $\nu\in \mathfrak a^*$, the class-one Whittaker function $\Psi_\nu$ (as a function 
on $\mathfrak a_0$) 
satisfies the partial differential equation
\begin{equation}\label{pde3}
 \frac{1}{2} \Delta f(x) - \sum_{\a\in\Pi} |\eta_\a|^2 e^{-2\alpha(x)} f(x) = \frac{1}{2} (\nu,\nu) f(x).
\end{equation}
For $\nu\in U$, this can be seen directly via the recursion~(\ref{recursion})
for the coefficients in the series expansion of the fundamental Whittaker functions $\Phi_\nu$.
In~\cite[Lemma 7.1]{h} it was shown that, for $\nu\in D$,
\begin{equation}\label{lim}
\lim_{x\to\infty,x\in\Omega} e^{s_0\nu(x)} \Psi_{\nu}(x) = c(\nu).
\end{equation}
By lemma~\ref{bounded}, if $\nu\in D$, then $e^{s_0\nu(x)} \Psi_{\nu}(x)$ is 
uniformly bounded for $x\in\mathfrak a_0 $.
Recalling Corollary~\ref{k}---note that the proof of uniqueness given there is valid for $\nu\in D$---we 
deduce the following characterisation of $\Psi_\nu$.
\begin{prop}\label{qtl}
For $\nu\in D$, the class-one Whittaker function $\Psi_\nu$ (on $\mathfrak a_0$) 
is the unique solution to (\ref{pde3})
such that $e^{s_0\nu(x)} \Psi_{\nu}(x)$ is bounded and (\ref{lim}) holds.
\end{prop}

\subsection{Weyl-invariant class-one Whittaker functions and an alternating sum formula}

In this section we present a variation of the formula (\ref{psi}) which generalises a formula given 
in~\cite{kl} for the case $G=SL(n,\R)$ and leads naturally to a normalisation for the class-one Whittaker 
functions which is invariant under the Weyl group $W$.  Using this, we also confirm a conjecture of 
Stade~\cite{stade} that a class-one Whittaker function can be expressed as an alternating sum of 
appropriately normalised fundamental Whittaker functions.  

Let 
\[a(\nu)=\prod_{\a\in\Sigma_+^\circ} \frac12 \left(|\eta_\a|/\sqrt{2(\a,\a)}\right)^{-\nu_\a}\G(\nu_\a).\]
\begin{prop}\label{succ} For $\nu\in U$,
\[c(\nu)^{-1} \Psi_\nu(x)= a(\nu)^{-1} \sum_{s\in W} a(s_0s\nu)\Phi_{s\nu}(x).\]
\end{prop}
\begin{proof}
From (\ref{psi}) we have
\[\Psi_\nu = \sum_{s\in W} M(s_0s,\nu)c(s_0s\nu)\Phi_{s\nu} .\]
It therefore suffices to show that, for all $s\in W$,
\[M(s,\nu) c(s\nu) a(s\nu)^{-1} = c(\nu) a(\nu)^{-1}.\]
We prove this by induction on $l(s)$.  If $s=s_\alpha,\ (\a\in\Pi)$, we have
\[M(s_\a,\nu)=\left(|\eta_\a|/2\sqrt{2(\a,\a)}\right)^{2\nu_\a} e_\alpha(\nu) e_\alpha(-\nu)^{-1},\]
\[c(s_\a\nu)=f_\a(-\nu)f_\a(\nu)^{-1} c(\nu) ,\]
\[a(s_\a\nu)^{-1}= \left(|\eta_\a|/\sqrt{2(\a,\a)}\right)^{-2\nu_\a}
\G(\nu_\a)\G(-\nu_\a)^{-1} a(\nu)^{-1}.\]
Using the duplication formula
\begin{equation}\label{dup}
\G(z)\G(z+\frac{1}{2})=2^{1-2z}\sqrt{\pi}\G(2z),
\end{equation}
we can write
\[e_\a(\nu)/f_\a(\nu)=\pi^{-1/2} 2^{\nu_\a-1+m(\a)/2} /\G(\nu_\a),\]
and so
\[M(s_\a,\nu)c(s_\a\nu)=
\left(|\eta_\a|/\sqrt{2(\a,\a)}\right)^{2\nu_\a} \frac{\G(-\nu_\a)}{\G(\nu_\a)} c(\nu).\]
Thus, 
\[M(s_\a,\nu)c(s_\a\nu)a(s_\a\nu)^{-1}=c(\nu)a(\nu)^{-1},\]
and the claim is proved for $l(s)=1$.  
For $s\in W$ and $\a\in\Pi$ with $l(s_\a s)=l(s)+1$,
\begin{eqnarray*}
\lefteqn{ M(s_\a s,\nu) c(s_\a s\nu) a(s_\a s\nu)^{-1} }\\
&=& M(s,\nu)M(s_\a ,s\nu)c(s_\a s\nu)a(s_\a s\nu)^{-1} \\
&=& M(s,\nu) c(s\nu) a(s\nu)^{-1} \\
&=& c(\nu) a(\nu)^{-1} ,
\end{eqnarray*}
by the induction hypothesis.
\end{proof}
Consider the normalised Whittaker functions
\[w_\nu(x) = a(\nu) c(\nu)^{-1}\Psi_\nu(x) \qquad \nu\in\mathfrak a^*,\ x\in\mathfrak a;\]
\[m_\nu(x)= \prod_{\a\in\Sigma_+^\circ}
\left(|\eta_\a|/\sqrt{2(\a,\a)}\right)^{\nu_\a} \G(1+\nu_\a)^{-1} \Phi_\nu(x) \qquad \nu\in U,\ x\in\mathfrak a.\]
By the above proposition and the functional equation
\begin{equation}\label{fun}
\G(z)\G(1-z)=\frac{\pi}{\sin \pi z} ,
\end{equation}
we have:
\begin{cor}\label{altsum}  For $\nu\in U$,
\[w_\nu(x) = R(\nu)^{-1} \sum_{s\in W} (-1)^{l(s_0s)} m_{s\nu}(x) ,\]
where
\[R(\nu)=\prod_{\a\in\Sigma_+^\circ} \frac{2 \sin \pi \nu_\a}{\pi} .\]
In particular, $w_\nu$ satisfies the functional equation
\[w_{s\nu}(x)=w_\nu(x),\qquad s\in W.\]
\end{cor}
This confirms a conjecture of Stade~\cite{stade}, who obtained this formula for the 
case $SL(3,\R)$ and conjectured that such a formula holds for all $SL(n,\R)$. 
In the case $G=SL(n,\R)$, the functions $w_\nu$ are 
essentially the same as those considered in~\cite{kl}.  

\subsection{The type $A_1$ case.}
Let $G=SL(2,\R)$.  Then we can identify $\mathfrak a_0$ with $\R$, and take
$\Sigma=\{\pm 1\}$, $\Pi=\{1\}$ and $m(1)=1$. Let $|\eta_1|^2=1/2$.
Then $L=2\Z_+$.  For $\l=2n$, write $c_n=c_\l(\nu)$. The recursion (\ref{recursion})
becomes $4(n^2+\nu n) c_n=c_{n-1}$ with $c_0=1$.  The solution is given by
\[c_n=\frac{4^{-n}\G(\nu+1)}{n!\G(n+\nu+1)},\]
and so
\[\Phi_\nu(x) = 2^{\nu} \G(1+\nu) \sum_{n\ge 0} \frac{(e^{-x}/2)^{2n+\nu}}{n!\G(n+\nu+1)} = 
2^{\nu} \G(1+\nu) I_\nu(e^{-x}),\]
where $I_\nu$ is the modified Bessel function of the first kind.
In this case, $W\simeq\Z_2$ acts on $\R$ by multiplication.  By the duplication
formula (\ref{dup}), we have
\[M(s_1,\nu)=4^{-\nu}\frac{\G(-\nu+\frac{1}{2})}{\G(-\nu+\frac{1}{2})},
\qquad c(\nu)=\frac{\sqrt{2\pi}\G(\nu)}{\G(\nu+\frac12)}.\]
Thus, using the functional equation~(\ref{fun}),
we obtain
\[\Psi_\nu(x) = M(s_1,\nu)c(-\nu)\Phi_\nu(x)+c(\nu)\Phi_{-\nu}(x) 
= \frac{2^{1-\nu}\sqrt{2\pi}}{\G(\nu+\frac12)} K_\nu(e^{-x}),\]
where
\[K_\nu(z) = \frac\pi 2 \frac{I_{-\nu}(z)-I_\nu(z)}{\sin\pi\nu}\]
is the Macdonald function.  Note that $a(\l)=2^{\l-1}\G(\l)$ and
the normalised Whittaker functions are given by 
$m_\nu(x)=I_\nu(e^{-x})$ and $w_\nu(x)=K_\nu(e^{-x})$.

\subsection{The type $A_2$ case.}\label{a2}
In this case we can identify $\mathfrak a_0$ with 
$\R^3_0=\{ x \in \mathbb{R}^3, x_1+x_2+x_3=0 \}$ and
take $\Pi=\{\a_1=(e_1-e_2)/\sqrt{2},\ \a_2=(e_2-e_3)/\sqrt{2}\}$, 
where $\{e_1,e_2,e_3\}$ is the standard basis 
for $\R^3$.  Set $m(\a_1)=m(\a_2)=1$,
$m(2\a_1)=m(2\a_2)=0$ and $|\eta_{\a_1}|^2=|\eta_{\a_2}|^2=2$.
For $\nu\in\R(\Pi)$ and $\l=2n\a_1+2m\a_2\in L=2\Z_+(\Pi)$,
write $c_{n,m}=c_\l(\nu)$.
Set $a=\a_1(\nu)$ and $b=\a_2(\nu)$.
Then the recursion (\ref{recursion}) becomes
\[(n^2+m^2-nm+an+bm) c_{n,m}=c_{n-1,m}+c_{n,m-1},\]
where $c_{0,0}=1$ and $c_{n,m}=0$ for $(n,m)\notin\Z_+^2$.
The solution is given by the following
formula, due to Bump~\cite{bump}:
\[c_{n,m}=\frac{\G(a+1)\G(b+1)\G(a+b+1)\Gamma(n+m+a+b+1)}
{n!m!\G(n+a+1)\G(m+b+1)\G(n+a+b+1)\G(m+a+b+1)}.
\]
In the notation of~\cite{bump,bh},
\[w_\nu(x)= \frac{\pi^2}{2} (y_1y_2)^{-1}W_{(\nu_1,\nu_2)}(y_1,y_2),\]
where 
\[\nu_1=(a+1)/3, \quad \nu_2=(b+1)/3,\quad y_1=2 e^{-\a_1(x)},\quad 
y_2= 2 e^{-\a_2(x)}.\]
The following integral representation is due to 
Vinogradov and Takhtadzhyan~\cite{vt}:
\begin{equation}\label{int}
w_\nu(x)= \frac12 (y_1/y_2)^{\frac{a-b}3} 
 \int_0^\infty K_{a+b}(y_1\sqrt{1+r}) K_{a+b}(y_2\sqrt{1+1/r} ) r^{\frac{a-b}2} \frac{dr}{r} .
\end{equation}
For $a=b=2/3$, we have the following simplification:
\[W_{(5/9,5/9)}(y_1,y_2)= 
\frac{2}{\sqrt{3}\pi} (y_1y_2)^{\frac13} (y_1^{\frac23}+y_2^{\frac23})^{\frac12}
K_{1/3}\left( (y_1^{\frac23}+y_2^{\frac23})^{\frac32} \right). \]
Using the integral representation (\ref{int}), Bump and Huntley~\cite{bh} derived an 
asymptotic expansion of $W_{(\nu_1,\nu_2)}(y_1,y_2)$ which is valid for large values 
of $y_1$ and $y_2$.  The leading term in the expansion is independent of the parameter $\nu$ 
and given by
\begin{equation}\label{expansion}
\sqrt\frac{2}{3\pi} (y_1y_2)^{\frac13} (y_1^{\frac23}+y_2^{\frac23})^{-\frac14}
\exp\left( - (y_1^{\frac23}+y_2^{\frac23})^{-\frac32} \right). 
\end{equation}
From this we deduce the following lemma, which we record for later reference.
\begin{lem}\label{asymp}
Let $\l_1,\l_2>0$.  If $y_1,y_2\to\infty$ with $y_2/y_1\to\delta$, then
\[\frac{ W_{(\nu_1,\nu_2)}(\sqrt{y_1^2+2\l_1y_1},\sqrt{y_2^2+2\l_2y_2}) }
{W_{(\nu_1,\nu_2)}(y_1,y_2)} \to \exp(-\l_1\varphi(\delta)-\l_2\varphi(1/\delta)),\]
where 
\[\varphi(d)=(1+d^\frac23)-d^\frac23(1+d^\frac23)^\frac12+d^\frac13 (1+d^{-\frac23})^\frac12.\]
\end{lem}

\subsection{Asymptotics for large $x$}

Consider the analytic function on $\mathfrak a^* \times \mathfrak a$ defined by 
\[\phi(\nu,x)=h(\nu)^{-1} \sum_{s\in W} (-1)^{l(s)} e^{s\nu(x)},\]
where $h(\nu)=\prod_{\a\in\Sigma_+^\circ} \nu_\a$.
Set \[\Omega^*=\Re(D)=\{\nu\in\mathfrak a_0^*:\ \nu_\a>0,\ \forall\a\in\Pi\}\]
and
\[\Omega=\{x\in\mathfrak a_0:\ \a(x)>0,\ \forall\a\in\Pi\}.\]
\begin{prop}\label{largex} 
Let $q=|\Sigma_+^\circ|$.
For all $x\in \Omega$ and $\nu\in\Omega^*$,
\[\lim_{c\downarrow 0} (2c)^q w_{-c\nu}(x/c)=\phi(\nu,x).\]
\end{prop}
\begin{proof}
First note that, since $\nu\in\Omega$, $cs\nu\in U$ for all $s\in W$ and
for all $c>0$ sufficiently small.
The claim follows from Corollary~\ref{altsum} and the fact (see~\cite{h})
that there exists a constant $k$ such that for all $s\in W$ and $c>0$ sufficiently small,
\[ \left| \sum_{\l\in L\backslash\{0\}} c_\l(-css_0\nu) e^{-\l(x)/c} \right| 
\le \sum_{n\ge 1} \frac{(n+d-1)!}{(d-1)!n!}\frac{k^n}{(n!)^2} e^{-2\min_{\a\in\Pi}\a(x)/c} .\]
\end{proof}

\section{Whittaker functions and exponential functionals of Brownian motion}

Define $k_\l=w_{-\l}$ for $\l\in\mathfrak a^*$.
Throughout this section we will 
identify $\mathfrak a_0^*$ with $\mathfrak a_0$ via the Killing form
and note that $\Omega^*=\Omega$.
Let $B^{(\mu)}$ be a Brownian motion in $\mathfrak a_0$ 
with covariance given by the Killing form and drift $\mu\in\Omega$.  
Then, by Corollary~\ref{k} and Proposition~\ref{qtl}, we have:
\begin{prop}\label{ew}
\[\mathbb E \exp\left( -\sum_{\a\in\Pi}  |\eta_\a|^2 e^{-2\alpha(x)} 
\int_0^\infty e^{-2\a(B^{(\mu)}_t)} dt \right) = e^{-\mu(x)} c(-s_0\mu)^{-1} \Psi_{-s_0\mu}(x) \]
\[=e^{-\mu(x)}  \prod_{\a\in\Sigma_+^\circ} 2 \left(|\eta_\a|/\sqrt{2(\a,\a)}\right)^{\mu_\a}
\G(\mu_\a)^{-1} k_{\mu}(x).\]
\end{prop}
In this context, the diffusion considered in section 3 has generator given by
\[\mathcal L_\mu = \frac12\Delta + \nabla\log k_{\mu} \cdot \nabla .\]
Note that this is well-defined for all $\mu\in\overline\Omega$.
Set
\[V_\mu(x) =  2 \sum_{\a\in\Pi} |\eta_\a|^2 e^{-2\alpha(x)} + (\mu,\mu) ,\]
and write $V=V_0$.  It follows from the intertwining
\begin{equation}\label{i1}
k_\mu\L_\mu=\frac{1}{2}(\Delta-V_\mu) k_\mu ,
\end{equation}
that the heat semigroup associated with $\L_\mu$ is given by
\begin{equation}\label{hk1}
P_t^\mu = k_\mu^{-1} Q_t^\mu k_\mu ,
\end{equation}
where $(Q_t^\mu)$ is the heat semigroup associated with $\frac12(\Delta-V_\mu)$.

Let $\mu\in\bar\Omega$ and consider the operator $\Lambda_\mu$, defined 
(on a suitable domain) by 
\[\Lambda_\mu e_\l = k_{\mu+\l},\qquad \l\in i\mathfrak a_0^*,\]
where $e_\l(x)=e^{\l(x)}$.  
Set $\K_\mu=k_{\mu}^{-1}\Lambda_\mu$.
In the type $A_1$ and $A_2$ cases, for each $x$, $k_{\mu+\l}(x)$ is a non-negative definite 
function of $\l$ and hence $\K_\mu$ is a Markov operator.  For the type $A_1$ case,
this follows from the integral representation~(\ref{mac}) and, for the type $A_2$ case,
it follows from the integral representation~(\ref{int}).  We remark that in fact it can be seen 
from an integral formula of Givental~\cite{givental} in the type $A_n$ case.
We conjecture that $\K_\mu$ is a Markov operator, in general.
In the type $A_1$ case, it is shown in~\cite{my} that $\K_\mu$ intertwines the
semigroup associated with $\L_\mu$ with the semigroup of a Brownian motion
with drift $\mu$.  This intertwining relation extends to the general setting:
\begin{prop} On a suitable domain, 
\[\L_\mu \K_\mu = \K_\mu (\frac{1}{2}\Delta+\mu\cdot\nabla) .\]
\end{prop}
\begin{proof}
For each $\l\in i\mathfrak a_0^*$, we have
\begin{align*}
(\Delta-V_\mu)\Lambda_\mu e_\l &= (\Delta-V_\mu)k_{\mu+\l} \\
& = (V_{\mu+\l}-V_\mu) k_{\mu+\l} \\
&= ( (\l,\l)+2(\mu,\l)) k_{\mu+\l} \\
&= ((\l,\l)+2(\mu,\l)) \Lambda_\mu e_\l \\
&= \Lambda_\mu(\Delta+2\mu\cdot\nabla) e_\l .
\end{align*}
Thus,
\[(\Delta-V_\mu)\Lambda_\mu=\Lambda_\mu(\Delta+2\mu\cdot\nabla).\]
Combining this with (\ref{i1}), we are done.
\end{proof}

\subsection{Brownian motion in a Weyl chamber and Duistermaat-Heckman measure}

Let $\mu\in\Omega$ and, for $c>0$, define $k^c_\mu(x)=k_{c\mu}(x/c)$.
By Proposition~\ref{largex}, the diffusion with generator
\[\L_\mu^c= \frac12\Delta + \nabla\log k^c_{\mu} \cdot \nabla ,\]
converges weakly as $c\downarrow 0$ to a Brownian motion with drift $\mu$ conditioned
(in the sense of Doob) never to exit the Weyl chamber $\Omega$ (see~\cite{bbo} for a
definition of this process).  In the limiting case $\mu=0$, the generator of the
Brownian moton conditioned never to exit $\Omega$ is given by
\[ \frac12\Delta + \nabla\log h \cdot \nabla ,\]
where $h(x)=\prod_{\a\in\Sigma^\circ_+}\a(x)$.
Note also that, as $c\downarrow 0$,
\[\Lambda^c e_\l (x) := (2c)^q \Lambda_0 e_{c\l} (x/c) = 
(2c)^q k_{c\l} (x/c) \to \phi(\l,x).\]
Thus, the intertwining operator $\Lambda^c$ converges, in a weak sense, to a positive integral 
operator with kernel given by $L(x,dt)=m^x_{DH}(dt)$, where $m^x_{DH}$ is the 
Duistermaat-Heckman measure associated with the point $x\in\Omega$, characterised by
\[\int_{\mathfrak a_0} e^{\l(t)} m^x_{DH}(dt) = \phi(\l,x),\qquad \l\in\mathfrak a^*. \]
This operator is discussed in~\cite{bbo}.  The intertwining
\[ (\Delta + 2 \nabla\log h \cdot \nabla ) L = L \Delta \]
plays a meaningful role in the multi-dimensional generalisations of Pitman's $2M-X$ theorem
obtained in~\cite{bbo,bj,oy}.  The operator $\Lambda$ has recently been shown to play
a similar role in the multi-dimensional version of the theorem of Matsumoto and Yor
obtained in~\cite{noc10} for the type $A_n$ case; it is a positive integral
operator with a kernel which can be interpreted as a kind of `tropical' analogue of the
Duistermaat-Heckman measure.

\section{The type $A_2$ case}

Consider the type $A_2$ case, as in section \ref{a2}.
For $x\in\R_0^3$, $\a_1(x)=(x^1-x^2)/\sqrt{2}$ and $\a_2(x)=(x^2-x^3)/\sqrt{2}$.
The Weyl chamber is $\Omega=\{x\in\R_0^3:\ x^1>x^2>x^3\}$.
Let $B^{(\mu)}$ be a Brownian motion in $\R_0^3$
with drift $\mu\in\Omega$. For $0\le t\le \infty$, set 
\[A^i_t=\int_0^t e^{-2\a_i(B^{(\mu)}_s)} ds ,\qquad i=1,2.\]
Let $\nu=-s_0\mu=(-\mu^3,-\mu^2,-\mu^1)$. Then, in the notation of section \ref{a2},
\[a=\frac{\nu^1-\nu^2}{\sqrt{2}},\quad b=\frac{\nu^2-\nu^3}{\sqrt{2}},\quad
\nu_1=\frac{a+1}3, \quad \nu_2=\frac{b+1}3.\]
Note that, for $x\in\R_0^3$,
\begin{align*}
\mu(x) &= \sqrt{2}(\nu^1+\nu^2)\a_1(x)+ \sqrt{2}\nu^1\a_2(x)\\
&= \frac{2a+4b}3\a_1(x) + \frac{4a+2b}3\a_2(x)\\
&= (2\nu_1+4\nu_2-2)\a_1(x)+ (4\nu_1+2\nu_2-2)\a_2(x).
\end{align*}
By Proposition~\ref{ew} and the integral formula~(\ref{int}),
\begin{eqnarray*}
\lefteqn{ \mathbb E \left(\exp\left( - \frac12 y^2_1 A^1_\infty - \frac12 y^2_2 A^2_\infty \right)\right) }\\
&=& 4\pi^2 y_1^{2\nu_1+4\nu_2-3} y_2^{4\nu_1+2\nu_2-3}
  \frac{2^{-a-b} }{\G(a)\G(b)\G(a+b)}
W_{(\nu_1,\nu_2)}(y_1,y_2) \\
&=& \frac{2^{2-a-b} }{\G(a)\G(b)\G(a+b)}
   \int_0^\infty y_1^{a+b-1} K_{a+b}(y_1\sqrt{1+r}) 
    y_2^{a+b-1} K_{a+b}(y_2\sqrt{1+1/r} ) r^{\frac{a-b}2} \frac{dr}{r} .
 \end{eqnarray*}
Let us observe that this Laplace transform can be inverted. Indeed, by using the fact that
\[
K_{a+b}(x)=\int_0^{+\infty} e^{-x \cosh u} \cosh ( (a+b)u) du
\]
we obtain
\begin{align*}
 \mathbb E \left(\exp\left( - \frac12 y^2_1 A^1_\infty - \frac12 y^2_2 A^2_\infty \right)\right)  
&=\frac{2^{2-a-b} }{\G(a)\G(b)\G(a+b)}(y_1y_2)^{a+b-1} \\
& \times
 \int_0^{+\infty}\int_0^{+\infty} \int_0^{+\infty}
 r^{\frac{a-b}2} e^{-y_1\cosh u \sqrt{1+r}} e^{-y_2\cosh v \sqrt{1+1/r}} \\
 & \times \cosh((a+b)u)\cosh((a+b)v)\frac{dr}{r} du dv
\end{align*}
and therefore, denoting $p$ the density of $(A_1^\infty, A_2^\infty)$, we get
\begin{align*}
p(y_1,y_2) =  & \frac{ 4}{\pi \G(a)\G(b)\G(a+b)} (2y_1y_2)^{-(a+b+1)/2} \\
&\times \int_0^{+\infty}\int_0^{+\infty} \int_0^{+\infty} r^{\frac{a-b}2} 
e^{-\frac{(1+r)\cosh^2u}{4y_1}-\frac{(1+1/r)\cosh^2 v}{4y_2}} \\
&  \times D_{a+b} \left( \frac{\sqrt{1+r}\cosh u}{\sqrt{y_1}}\right)
 D_{a+b} \left( \frac{\sqrt{1+1/r}\cosh v}{\sqrt{y_2}}\right) \\
 & \times \cosh((a+b)u)\cosh((a+b)v)\frac{dr}{r} du dv .
\end{align*}

\subsection{The intertwining kernel}

Suppose $\mu=0$.  The intertwining operator $\Lambda=\Lambda_0$ satisfes
$\Lambda e_{-\nu} = w_\nu$.
Let $a=\a_1(\nu)$, $b=\a_2(\nu)$ and write $t=t_1\a_1+t_2\a_2$.
By the integral formula (\ref{int}), 
\begin{eqnarray*}
\lefteqn{ w_\nu(x) = \frac12 \int_0^\infty\int_0^\infty\int_0^\infty
 \left( \delta^{-\frac13} \frac{\sqrt{r}}{uv}\right)^a
\left( \delta^{\frac13} \frac{1}{uv\sqrt{r}}\right)^b }\\
&& \times \exp\left( -\frac{y_1\sqrt{1+r}}{2}(u+\frac1u ) -\frac{y_2\sqrt{1+1/r}}{2}(v+\frac1v )\right)
\frac{du}u\frac{dv}v\frac{dr}r ,
\end{eqnarray*}
where $y_1=2 e^{-\a_1(x)}$, $y_2= 2 e^{-\a_2(x)}$ and $\delta=y_2/y_1$.  
It follows, by a straightforward calculation, that we can write
\[\Lambda e_{-\nu} (x) = \int_{-\infty}^\infty\int_{-\infty}^\infty
 \Lambda(x,t) e^{-at_1} e^{-bt_2} dt_1 dt_2 ,\]
where
\[\Lambda(x,t)= K_0\left( \sqrt{
y_1^2 (1+\delta^{\frac23}e^{-t_1+t_2}) (1+\delta^{\frac23}e^{t_1}) (1+\delta^{\frac23}e^{-t_2})}
\right)
.\]
A plot of $\Lambda(x,t)$, with $x$ fixed, is shown in Figure 1.  As explained in~\cite{noc10},
the measure with density given by $\Lambda(x,\cdot)$ can be interpreted as a kind of `tropical' analogue of the Duistermaat-Heckman measure associated with the point $x\in\Omega$.

\begin{figure}
\centering
\includegraphics[width=10cm]{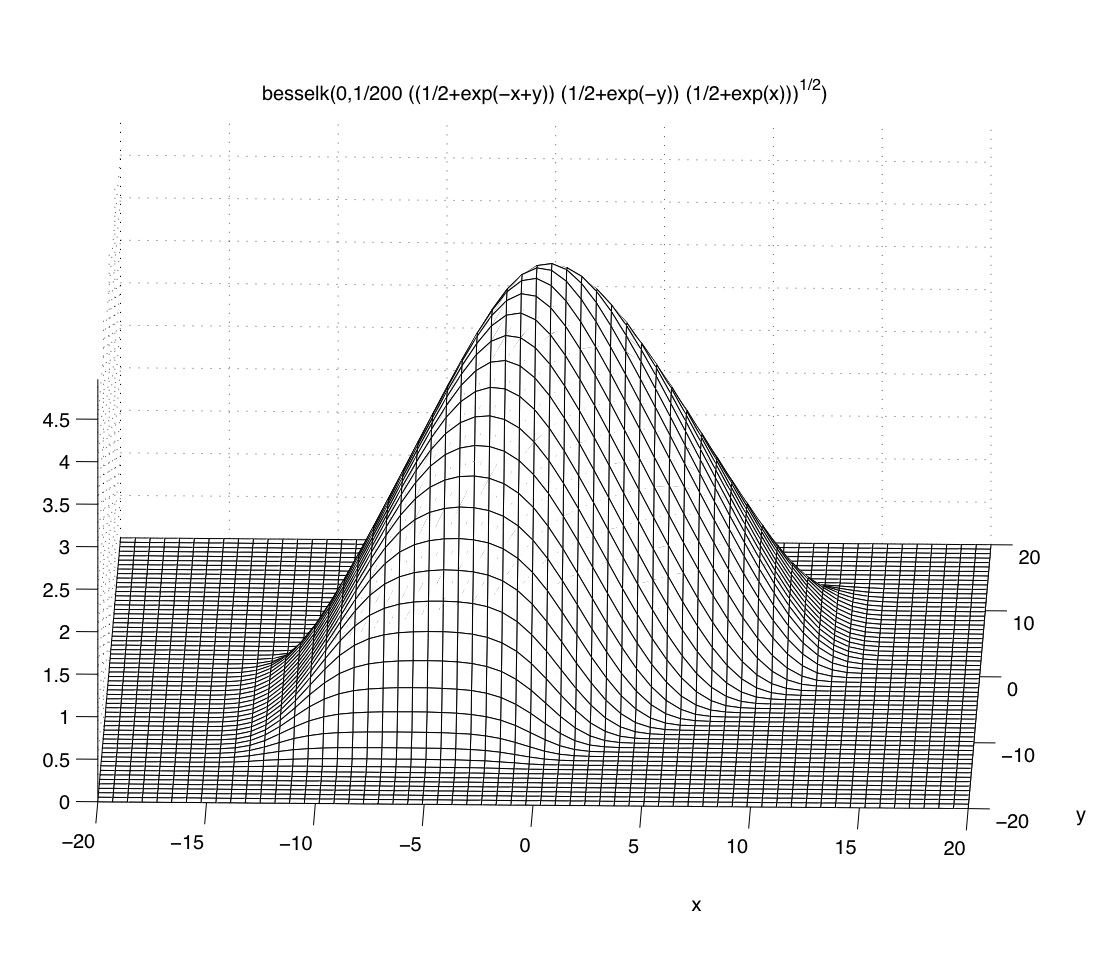}
\caption{The intertwining kernel $\Lambda(x,\cdot)$}
\end{figure}

\subsection{Behaviour at $-\infty$}

From the asymptotic expansion of~\cite{bh}, for any $\l$ we have 
$k_\mu(x)^{-1} k_{\lambda}(x)\to 1$
as $x\to-\infty$ (in the sense that $\alpha_1(x)\to-\infty$ and $\alpha_2(x) \to-\infty$).
This suggests that the process with generator
\[\mathcal L_\mu = \frac12\Delta + \nabla\log k_\mu \cdot \nabla \]
has a unique entrance law starting from $-\infty$, given by 
\[p^\mu_t(dx) = e^{-\frac12\|\mu\|^2 t} k_\mu(x) \theta_t(dx),\]
where, for each $t>0$, 
\[\int k_{i\tau}(x) \theta_t(dx) = e^{-\frac12\|\tau\|^2 t},\qquad \tau\in\Omega .\]
The existence of this entrance law is established (more generally, for type $A_n$)
in the paper~\cite{noc10}.  It is clear from the proof given there that it should be unique.
On the other hand:
\begin{proposition}
Let $(X_t^{x_0})_{t \ge 0}$ be the diffusion with generator $\mathcal L_\mu$ started at $x_0$. 
If $\alpha_1(x_0)\to-\infty$ and $\alpha_2(x_0) \to-\infty$, with  
$\alpha_1 (x_0)-\alpha_2 (x_0) \rightarrow \kappa$, 
then
\[
\left(e^{\alpha_1(x_0)}\int_0^{+\infty} e^{-2\alpha_1(X_s^{x_0})}ds,  e^{\alpha_2(x_0)}\int_0^{+\infty} e^{-2\alpha_2(X_s^{x_0})}ds \right)
\]
converges in probability to $\left( \varphi (e^\kappa), \varphi (e^{-\kappa})\right)$, where
\[\varphi(d)=(1+d^\frac23)-d^\frac23(1+d^\frac23)^\frac12+d^\frac13 (1+d^{-\frac23})^\frac12.\]
\end{proposition}
\begin{proof}
Let $\l_1,\l_2>0$.  
We easily compute
\begin{align*}
 & \mathbb{E} \left( \exp \left(-\l_1e^{\alpha_1(x_0)}\int_0^{+\infty} e^{-2\alpha_1(X_s^{x_0})}ds-\l_2  e^{\alpha_2(x_0)}\int_0^{+\infty} e^{-2\alpha_2(X_s^{x_0})}ds \right) \right) \\
 = &\frac{ \mathbb E \exp\left( - \frac12 (y_1^2+2\l_1y_1) A^1_\infty - \frac12 
(y_2^2+2\l_2y_2) A^2_\infty \right) }
{\mathbb E \exp\left( - \frac12 y^2_1 A^1_\infty - \frac12 y^2_2 A^2_\infty \right)}
\end{align*}
with $y_1=2e^{-\alpha_1 (x_0)}, y_2=2e^{-\alpha_2 (x_0)}$.
But, by Lemma~\ref{asymp}, if $y_1,y_2\to\infty$ with 
$y_2/y_1\to\delta=e^{\kappa}$, then
\[\frac{ \mathbb E \exp\left( - \frac12 (y_1^2+2\l_1y_1) A^1_\infty - \frac12 
(y_2^2+2\l_2y_2) A^2_\infty \right) }
{\mathbb E \exp\left( - \frac12 y^2_1 A^1_\infty - \frac12 y^2_2 A^2_\infty \right)}
 \to e^{-\l_1\varphi(\delta)-\l_2\varphi(1/\delta)}.\]
\end{proof}

\bigskip

\noindent {\em Acknowledgements.} 
We are grateful to the anonymous referees for helpful suggestions.  
Research of second author was supported in part
by Science Foundation Ireland grant number SFI 04-RP1-I512.

\end{document}